\documentclass[12pt,a4paper,twoside,final,notitlepage, reqno]{article}

\usepackage[english]{babel}
\usepackage[latin1]{inputenc}
\usepackage[a4paper,left=2.5cm,right=2.5cm]{geometry}
\usepackage{amssymb}
\usepackage{amsthm} 
\usepackage{amsmath}
\usepackage{graphics,graphicx}
\usepackage[usenames, dvipsnames]{color}
\usepackage{booktabs}
\usepackage{subfigure} 
\usepackage{siunitx}

\graphicspath{{figs/}}

\makeatletter
\def\input@path{{figs/}}
\makeatother

\theoremstyle{definition}

\theoremstyle{definition}

\theoremstyle{definition}

\providecommand{\keywords}[1]  {\textbf{Keywords:} #1}

\providecommand{\acknow}[1] {\textbf{Acknowledgments.} #1}

\def\fs{ \footnotesize}
\def\ss{ \scriptsize}

\usepackage{esdiff}
\renewcommand{\d}[1]{\mathrm{d} #1}

\newcommand{\R}{\mathbb{R}}
\renewcommand{\phi}{\varphi}

\DeclareMathOperator{\e}{e}

\usepackage{authblk}

\title{
  \vspace{-1.5cm}
  \bf{\Large{
      Efficient high order schemes for stiff ODEs
      \\[5pt]
      in cardiac electrophysiology
    }}
}

\author[1]{Charlie Douanla-Lontsi \thanks{charlie.douanla-lontsi@inria.fr}}
\author[1,2]{Yves Coudi\`ere \thanks{yves.coudiere@inria.fr}}
\author[3]{Charles Pierre \thanks{charles.pierre@univ-pau.fr}}

\affil[1]{
  INRIA Bordeaux Sud Ouest, Universit\'e de Bordeaux, France.
}
\affil[2]{
  Institut de Math\'ematiques de Bordeaux, 
  UMR CNRS 5241.
}\affil[3]{
  Laboratoire  de Math\'ematiques et de leurs Applications,
  UMR CNRS 5142, \protect \\
  Universit\'e de Pau et des Pays de l'Adour, France.
}

\usepackage{fancyhdr}

\begin{document} 

\date{13 October, 2017}

\maketitle

\begin{abstract}
  In this work we analyze the resort to high order exponential solvers for stiff ODEs in the context of cardiac electrophysiology modeling. The exponential Adams-Bashforth and the Rush-Larsen schemes will be considered up to order 4. These methods are explicit multistep schemes.The accuracy and the cost of these methods are numerically analyzed in this paper and benchmarked with  several classical explicit and implicit schemes
  at various orders. This analysis has been led considering data of high particular interest in cardiac electrophysiology : the activation time ($t_a$ ), the recovery time ($t_r $) and the action potential duration ($APD$).
  The Beeler Reuter ionic model, especially designed for cardiac ventricular cells, has been used for this study.
  It is shown that, in spite of  the stiffness of the considered model, exponential solvers allow computation at large time steps, as large as for implicit methods. Moreover, in terms of cost for a given accuracy, a significant gain is achieved with exponential solvers.
  We conclude that accurate computations at large time step are possible with explicit high order methods. This is a quite important feature when considering stiff non linear ODEs.
\end{abstract}
\vspace{20pt}
\noindent
  \\[3pt]
  \acknow{
  This study received financial support from the french government as part of the
  ``Investissement d'Avenir'' program managed by the ``Agence Nationale de la Recherche''
  (ANR), grant reference ANR-10-IAHU-04. It also received fundings of the  ANR project HR-CEM n0.\ 13-MONU-0004-04.
  \\ \\
}
\section{Introduction}
The numerical resolution of stiff ordinary differential equations (ODEs) is an issue encountered in many fields of applied sciences. In cardiac electrophysiology, the electrical activity of the heart is described by a system of parabolic partial differential equations coupled with a system of ODEs called \textit{ionic models}. 
The stiffness and the nonlinearity of the ionic models (see \cite{spiteri} for the stiffness analysis) make their numerical resolution very challenging. The classical schemes have serious drawbacks to solve such ODEs. On the one hand, the classical stable methods are implicit and lead to high computational cost (because of the nonlinear solvers) for large time-steps, on the other hand explicit solvers require very small time steps also leading to high computational costs. Meanwhile current solvers in cardiac electrophysiology are usually based on order 1 or 2 schemes (see \cite{perego-2009, GRL, RL1,ethier2008semi}).
In this paper we investigate the resort to a class of both explicit and stable schemes referred as «exponential methods » of high order as an alternative to solve cardiac electrophysiological problems. Namely we will consider the exponential Adams-Bashforth (EAB) and the Rush-Larsen (RL) techniques.\\
Let us consider the general initial value problem,
\begin{align}
  \label{F1}
  \diff{y}{t} = F(t,y) \ \ t \in (0,T], \quad 
  y(0) = y_{0} \in \R^N. 
\end{align}
EAB schemes \cite{EAB_M_O} and RL schemes \cite{coudiere-lontsi-pierre-RLk} take advantage of a splitting of the model function F into some linear part $a$ and a nonlinear part $b$, such that \eqref{F1} becomes,
\begin{align}
  \label{F2}
  \diff{y}{t} = a(t,y)y + b(t,y), \quad   y(0) = y_{0}\in \R^N. 
\end{align}
Notice that in \eqref{F2}, $a$ is not the exact linear part of $F$ (its differential) but, an approximation or a guess thereof. 
The EAB and RL are built from a transformation of \eqref{F2} on each time discretization interval $[t_n, t_{n+1}]$ in the following form,
\begin{align}
  \label{F3}
  \diff{y}{t} = \alpha_n y + c_n(t,y), \quad   y(0) = y_{0}\in \R^N. 
\end{align} 
Where $\alpha_n \in \R^N$ is a stabilizer set at every time step and $c_n(t,y) = (a(t,y) - a_n)y + b(t,y)$ . With the formulation \eqref{F3}, the exact solution satisfies the variation of the constant formula,
\begin{align}
  \label{var_of_const_for}
  y(t_{n+1})=\e^{\alpha_n h}\left( y(t_n)+ \int_{t_n}^{t_{n+1}} \e^{-\alpha_n(\tau - t_n)} c_n(\tau, y(\tau)) \d{\tau} \right).
\end{align}
EAB and RL schemes are based on this formula and are obtained by replacing $c_n$ by an approximation. Their precise definitions will be given in Section \ref{sch:statement}. 
The aim of this paper is to study the efficiency of EAB and RL methods of order 1 up to 4, and to show for a given scheme of order $k$, how to compute without degrading the accuracy, the activation time ($t_a$), the recovery time ($t_r$) and the action potential duration ($APD$) which are informative values of a particular interest in cardiac electrophysiology.
The efficiency of the schemes is analyzed both in terms of the accuracy and of the cost. The comparison is made using a realistic test case and is completed by including a benchmark with several classical methods either of implicit or explicit type: The Crank-Nicolson (CN), the Runge Kutta (RK$_4$), the Adams-Bashforth (AB$_k$), and the backward differentiation (BDF$_k$) (see \cite{ethier2008semi, Hairer-ODE-I} for the statement of these schemes). One interesting result from this test case is the possibility with EAB and RL, to use large time steps as in the case of implicit schemes but with the same cost as for the explicit schemes. 
The paper is organized as follows. In Section \ref{sch:statement} are presented the stabilized schemes. A brief description of the transmembrane action potential, $t_a$, $t_r$, $APD$ and ionic model is given in Section \ref{sc:Medeling}. The methodology used to compare the numerical schemes and to compute the values $t_a$, $t_r$ and $APD$ predicted by a numerical solution are developed in Section \ref{sc:tools}. The comparison of the methods follows in Section \ref{sc:num_res}. 
In Section \ref{sc:conclusion}, a brief conclusion is given.
\section{EAB$_k$ and RL$_k$ scheme statements}
\label{sch:statement}
When the function $c_n(t,y)$ in \eqref{F3} is a polynomial $\sum_{j=0}^{k-1} p_j(t-t_n)^j $ of degree $k-1$, the relation \eqref{var_of_const_for} becomes,
$$y_{n+1} = \e^{\alpha_n h} y_n + h \sum_{j=0}^{k-1}p_j j!h^j \varphi_{j+1}(\alpha_n h),$$
where the functions $\varphi_j$ are defined recursively by,
\begin{align*}
  \varphi_{j+1}(z) = \frac{\varphi_j(z) - \varphi_j(0)}{z}, \quad \varphi_{0}(z) =\e^z  \quad \text{and} \quad  \varphi_j(0)  = \frac1{j!} \quad \forall j \geq 0.
\end{align*}
The schemes introduced in the sequel are multi-steps. We will use the following notation $a_n=a(t_n, y_n)$, $b_n= b(t_n, y_n)$. 
\begin{table}[htbp!]
  \centering
  \caption{Coefficients $\gamma_{nj}$ for the EAB$_k$ schemes.}
  \label{tab:betanj}
  \begin{tabular}{lllll}\toprule
    $k$ & 1 & 2 & 3 & 4 \\ \midrule
    $\gamma_{n0}$ & $c_n^n$ & $c_n^n$ & $c_n^n$ & $c_n^n$ \\
    $\gamma_{n1}$ & -- & $c_n^n - c_n^{n-1}$ & $\frac{3}{2} c_n^n -2 c_n^{n-1} + \frac{1}{2} c_n^{n-2}$ & $\frac{11}{6}c_n^n-3c_n^{n-1}+\frac{3}{2}c_n^{n-2}-\frac{1}{3}c_n^{n-3}$ \\
    $\gamma_{n2}$ & -- & -- & $c_n^n-2c_n^{n-1} + c_n^{n-2}$ & $2c_n^n-5c_n^{n-1}+4c_n^{n-2}-c_n^{n-3}$ \\
    $\gamma_{n3}$ & -- & -- & -- & $c_n^n-3c_n^{n-1}+3c_n^{n-2}-c_n^{n-3}$ \\ \bottomrule
  \end{tabular}
\end{table}
\begin{itemize}
\item EAB$_k$: On the one hand we set $\alpha_n = a_n$, on the other hand the function $c_n$ in \eqref{F3} is approximated by its Lagrange  interpolation polynomial $\tilde{c}_n$ of degree $k-1$ at the time instants $t_n,\dots,t_{n-k+1}$. This polynomial satisfies $\tilde{c}_n(t_{n-j}) =c_n(t_{n-j},y_{n-j})$ for $j=0,\dots,k-1$.  The values $c_n(t_{n-j},y_{n-j})$ are given by $c_n^{n-j}= b_{n-j}+(a_{n-j}-a_n)y_{n-j}$ for $j=0,\dots,k-1$. If we write $\tilde{c}_n(t) = \sum_{j=0}^{k-1} \frac{\gamma_{nj}}{j!} \left(\frac{t-t_n}{h}\right)^j$, the definition of the EAB$_k$ scheme is deduced from the formula \eqref{var_of_const_for} by 
  \begin{equation}
    \label{eq:EAB-scheme-2}
    y_{n+1} = \e^{a_n h} y_n + h \sum_{j=0}^{k-1} \gamma_{nj}
    \varphi_{j+1}(a_n h),
  \end{equation} 
  where the coefficients $\gamma_{nj}$ are given in Table \ref{tab:betanj}.
\item RL$_k$: In the case the function $c_n(t,y)$ in \eqref{var_of_const_for} is a constant $c_n =\beta_n \in \R$ then we have the following simple scheme definition,
  \begin{align}
    \label{gen_RLk}
    y_{n+1} = y_n + h \varphi_1(\alpha_n h)(\alpha_n y_n + \beta_n),
  \end{align}
  that we refer as Rush-Larsen schemes as in the continuity of \cite{perego-2009}. The following choices for defining $\alpha_n$ and $\beta_n$ ensure the convergence at order k of the scheme \eqref{gen_RLk} and thus are named Rush-Larsen schemes of order k (RL$_k$). 

  \begin{itemize} 
  \item $k=1:\quad \alpha_n = a_n , \quad \beta_n = b_n$.
  \item $k=2:\quad \alpha_n = \frac{3}{2}a_n-\frac{1}{2}a_{n-1}, 
    \quad \beta_n  = \frac{3}{2}b_n-\frac{1}{2}b_{n-1}$.
  \item $k=3:\quad  \alpha_n = \frac{1}{12}(23a_n - 16 a_{n-1} + 5 a_{n-2})$,          
    \\[5pt]
    $\beta_n = \frac{1}{12}(23b_n - 16 b_{n-1} + 5 b_{n-2}) + \frac{h}{12}(a_nb_{n-1} - a_{n-1}b_n).$ 
  \item $k=4: \quad 
    \alpha_n = \frac{1}{24}(55a_n - 59 a_{n-1} + 37 a_{n-2}-9a_{n-3})$,         
    \\[5pt]
    $\beta_n = \frac{1}{24}(55b_n - 59 b_{n-1} + 37 b_{n-2}-9b_{n-3}) + \frac{h}{12}(a_n(3b_{n-1}-b_{n-2}) - (3a_{n-1}-a_{n-2})b_n).$
  \end{itemize}

\end{itemize}
Notice that the EAB$_1$ scheme is the same with RL$_1$ scheme and also the exponential Euler scheme.

The previous description of the EAB$_k$ scheme has been given very briefly but, more details can be found in \cite{EAB_M_O,NorsEAB} (for general ODEs) and in \cite{coudiere-lontsi-pierre-2016} for cardiac electrophysiology application. 
\section{Modeling in cellular cardiac electrophysiology}
\label{sc:Medeling}
\subsection{The action potential}
\label{sec:action-potential}
\begin{figure}[htbp!]
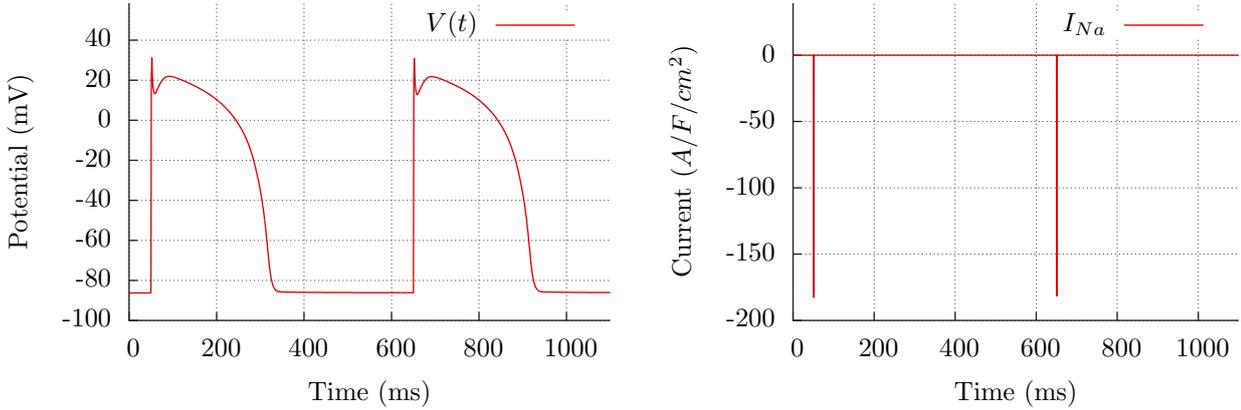

  \centering
  \input{V_tnnp.txt}~~~~~~~~\input{INa_tnnp.txt}\\[5pt]
  \caption{TNNP model \cite{tnnp} illustration. Left, two cellular action
    potentials : starting at a negative resting value, the
    transmembrane voltage $V(t)$ has a stiff depolarization followed
    by a plateau and repolarizing to the resting value. Right : each
    depolarization is induced by an ionic sodium current $I_{Na}(t)$}
  \label{fig:action-potential}
\end{figure}
The phenomenon studied here is the so called \textit{cellular action potential}, that we briefly present here. A potential difference is observed between the inside and outside of the cell, said membrane potential and denoted $V$. This potential caused by the differences in ionic concentrations between the inside and outside of the cells is dynamic in time, as well as these ionic concentrations. The potential $V$ can abruptly switch from a \textit{resting} state (during which $V=V_r\simeq -100 mV$ ) to an \textit{excited} state (where $V$ is in the range of 10 mV) in which it is maintained during a few tenth of seconds before returning to its resting state (see Figure \ref{fig:action-potential}).
It is this cycle,
\begin{displaymath}
  \text{resting state} 
  \overset{\text{excitation}}{\quad\quad \longrightarrow\quad\quad} 
  \text{excited state} 
  \overset{\text{recovery}}{\quad\quad\longrightarrow\quad\quad} 
  \text{resting state},
\end{displaymath}
that one designates as action potential. The resting potential $V_r$ is associated  to a pic potential $V_p$ corresponding to the maximum of the potential $V$ at the end of the excitation and a threshold potential $V_{th}$ such that $V_r < V_{th}<V_p$. 

We adopt here the following definitions : $V_{th}$ is the potential corresponding to $20\%$ of depolarization, the activation time $(t_a)$ and the recovery time $(t_r)$ are the time instants where the potential reaches the value $V_{th}$ the first and the second time respectively and the action potential duration $(APD)$ is the amount of time in which the voltage remains elevated above $V_{th}$. More precisely,
\begin{equation}
  \label{eq:def-APD}
  V_{th} = 0.8 V_r + 0.2 V_p, \quad V(t_a) = V_{th} = V(t_r), \ t_a < t_r \quad \text{and} \quad
  \text{APD} = t_r - t_a.
\end{equation}
\subsection{Ionic Models} 
\label{sec:mod-ion}
The variations of the ionic concentrations are described by ionic models and are systems of ODE. Experimental models (such as BR \cite{beeler-reuter} and TNNP \cite{tnnp} models designed for cardiac human cells) involve a variable $y\in \R^N$ composed of the following entries:
\begin{itemize}
\item \textbf{The membrane potential:} $V$ in mV. The equation on the potential is written,
  \begin{equation}
    \label{eq:edo-V}
    \diff{V}{t} = -I_{ion}(y(t)) + I_{st}(t),
  \end{equation}
  where $I_{ion}$ (reaction term) is the total ionic current crossing the membrane cell and $I_{st}$ is the stimulation current, it is a source term.
\item \textbf{The gating variables:} they are parameters between 0 and 1 expressing the variability and the permeability of the membrane cell for the specific ionic species. One denote by $W \in \R^P$ the vector of gating variables. The equations on $W$ are, for $i=1\dots P$,   
  \begin{equation}
    \label{eq:edo-gating-ver}
    \diff{W_i}{t} = \dfrac{W_{\infty,i}(y) - W_i}{\tau_i(y)},
  \end{equation}
  where $W_{\infty,i}(y)\in \R$, $\tau_i(y)\in\R$ are scalar functions given by the model. In these equations the linear and nonlinear parts are encoded in the model and are equal to $-1/\tau_i(y)$ and $W_{\infty,i}(y) / \tau_i(y)$ respectively.
\item \textbf{Ionic concentrations:} One denote by $C\in \R^{N-P-1}$ the vector of concentrations.
\end{itemize}
All these entries are collected in the vector $y$ as follows ,
\begin{equation*}
  y=\left [
    \begin{array}{c}
      W \\ X
    \end{array}
  \right. 
  , \quad 
  X=\left [
    \begin{array}{c}
      C \\ V
    \end{array}
  \right.
  , \quad W\in\R^P,\quad C\in\R^{N-P-1}, \quad  V=y_N\in\R, 
\end{equation*}
The sub-vectors $W$ corresponds to the lines of (\ref{F1}) including stabilization with the linear part $a(t,y) = -1/\tau(y)$ and the non linear part $b(t,y) = -W_\infty (y)/\tau(y)$. The sub-vector $X$ corresponds to the lines of (\ref{F1}) with no stabilization ( $a(t,y) = 0$). 
The associated ODE written in the form (\ref{F2}) is then defined by,
\begin{displaymath}
  a(t,y) = \left[
    \begin{array}{cc}
      A_1(t,y) & 0
      \\
      0 & 0
    \end{array}\right], 
  \quad b(t,y) = \left[
    \begin{array}{c}
      B_1(t,y) 
      \\
      B_2(t,y) 
    \end{array}\right.,
\end{displaymath}
where the matrix $A_1(t,y)\in\R^P\times \R^P$ is diagonal, $A_1(t,y)= \text{Diag}(-1/\tau_i(y))$, and $B_1(t,y) = \{W_{\infty, i}(y)/\tau_i(y), i=1\dots P\}\in\R^P$.

\section{Numerical Study}
\label{sc:tools}
\subsection{Scheme analysis methods}
\begin{itemize}
\item \textbf{Test case :} The evaluation and comparisons between different ODE solvers is done with a test case. Specifically, the Beeler Reuter \cite{beeler-reuter} model is considered and written in the form \eqref{F2} as described in Section \ref{sec:mod-ion}. We denote by $y(t)$ the solution of the associated ODE \eqref{F2} in $(0, T]$ with $T = 396 \ ms $. this solution is uniquely defined once the initial condition $y_0$ and the stimulation current $I_{st}$ in \eqref{eq:edo-V} are fixed. $y_0$ is the resting state as described by the model. 
  The function $I_{st}(t)$ is positive, null outside the interval $(t_s-1,t_s+1)$, $t_s$=20 ms and with integral $\int_0^T I_{st}(t) dt = I_{stim}$, a typical current of stimulation fixed by the models, equal to 50 mA. We also impose to $I_{st}$ a $C^4$ regularity in order to observe the convergence orders of schemes up to 4.
  
\item \textbf{Numerical solution:} Let $m\ge 1$ be an integer for which one associates the time-step $ h =T/m$ and the regular mesh $ T _m = \{t_n=j h ,~j=0\dots m\}$ of the interval $(0,T]$. The numerical solution $(y^n)$ is the element of the space $E_m$, $E_m=\{(y^n)_{0\le n\le m}, \ y^n\in\R^N\}$.
  The space $E_m$ of the numerical solutions is simply $(\R^N)^m$ but to $(y^n)\in E_m$ is implicitly associated a time-step $h$ and a mesh $ T _m$, such that each value $y^n$, $0\le n \le m$ of $(y^n)\in E_m$ is supposed to be an approximation of $y(t_n)$.

\item \textbf{Reference solution:} For a given test case, we cannot access to the exact solution $y(t)$ of the associated ODE. So for a numerical solution $(y^n)\in E_m$, we set $m' = 2^r m$ with $r\ge 0$ an integer and define the reference solution associated to $(y^n)$ (or $m$) as the numerical solution $y_{ref}\in E_{m'}$ for the problem \eqref{F1}, computed by the RK$_4$ scheme with the time-step $ h _{ref} = T/m' =  h  / 2^r$.
  The reference solution $y_{ref}$ is then not unique and depend on $r$. In practice $r$ is chosen \textit{large enough} such that the error between the exact solution $y$ and $y_{ref}$ is negligible compared to the error between the numerical solution $(y^n)$ and $y_{ref}$.

\item \textbf{Interpolation of the solution:}
  To compare the numerical solution with the reference solution and to be able to compute the numerical error in terms of function norm, we define an interpolator  $\pi_{m,i}: E_m \longrightarrow C^0(0,T],$ transforming the component $i$ of the numerical solution $(y^n)\in E_m$ in $C^0(0,T]$, the set of the continuous functions on $(0,T]$. We impose to the interpolant $\pi_{m,i} y^n$ to be a polynomial piecewise function of degree 3, this constraint is necessary to observe the convergence order up to 4. We assume that $m$ is a multiple of 3 and fix $(y^n)\in E_m$. We decompose the interval $[0,T]$ in a sequence of 3 intervals packages $P_s=[t_{3s},t_{3s+1}] \cup [t_{3s+1},t_{3s+2}] \cup [t_{3s+2},t_{3(s+1)}]$ , for $s=0\dots m/3$. The interpolated $f := \pi_{m,i} y^n$ is the unique polynomial of degree 3 on each $P_s$, continuous on $[0,T]$, such that $f(t_n) = y^n_i$ for all $n=0\dots m$. This interpolator is not Canonical: an $H^3$-Hermite interpolation on each interval $(t_n,t_{n+1})$ is an alternative.  
  The emphasis will be on the membrane potential $V(t)=y_N(t)$ and for more simplicity we note $\pi_m = \pi_{m,N}$ and $\pi = \pi_{m,N}$  in confusion absence. 
\item \textbf{Accuracy:} Let $(y^n)$ be a numerical solution and $y_{ref}$ a reference solution.
  We denote $\widehat{V}_{ref}$ and $\pi y^n = \widehat{V}$ the reference membrane potential and the membrane potential interpolating associated to $(y^n)$. The accuracy of each method is evaluated through a relative error between the reference solution and the numerical solution. 
  We define the errors in norm $ L^\infty $ by :
  \begin{equation}
    \label{eq:e_Lp-def}
    e_\infty  = \dfrac{\max \vert \widehat{V} - \widehat{V}_{ref} \vert}{\max \vert \widehat{V}_{ref} \vert}.
  \end{equation}
  Notice that the choice of the membrane potential $V$ is arbitrary and that any other component of $(y^n)$ could have been considered.
  The accuracy notion will be central here and it is convenient to identify several aspects.
\item \textbf{Cost:} The accuracy takes all its meaning when it is associated with a cost. Here it is a \textit{computational} cost and is evaluated with the CPU time during a simulation. It is evaluated by our fortran 90 code for each simulation. The CPU times depend on the computer used to perform the solutions and how the numerical solver is implemented. This is especially the case of implicit solvers where the Newton-Krylov algorithm type is used with the approximation of the Jacobian by the finite difference method.
\end{itemize}

\subsection{Computation of $t_a$, $t_r$, $APD$ and the associated errors}
\label{sc:phys_quantities}
Let $(y^n)$ be a numerical solution given and $V^n=y^n_N$, $0\le n\le m$, the membrane potential associated to $(y^n)$. Then there exist two unique indexes $n_a < n_r$ such that
\begin{displaymath}
  V_{n_a+1} \le V_{th} < V_{n_a+2} \quad ,\qquad 
  V_{n_r+1} \le V_{th} < V_{n_r+2},
\end{displaymath}
for the threshold potential $V_{th}$ defined in Section \ref{sec:action-potential}.

On the intervals 
$(t_{n_i},t_{n_i+3})$
, $i\in\{a,r\}$,  we compute the Lagrangian interpolation polynomial of degree 3  $p_i(t)$ for the values $V_j$ associated to $t_j$, $j=n_i,\dots,n_i +3$.  The activation time $t_a$ an the recovery time $t_r$ are the computed as the solution of,
\begin{displaymath}
  p_a(t_a) = V_{th}, \quad p_r(t_r) = V_{th}.
\end{displaymath}

Again the use of interpolation of order 3 is necessary to observe the convergence order up to 4. In above, we suppose that all is well defined, which is the case if the numerical solution $(y^n)$ is physiologically relevant. 

The relative error between the activation, recovery time and APD predicted by a numerical solution $(y^n)$ and a reference solution $y_{ref}$ will be computed by,
\begin{displaymath}
  e_{t_a} = \dfrac{ | t_a - t_{a,ref}|}{|t_{a,ref}|},\quad 
  e_{t_r} = \dfrac{|t_r - t_{r,ref}|}{|t_{r,ref}|}, \quad 
  e_{\text{APD}} = \dfrac{|\text{APD} - \text{APD}_{ref}|}{|\text{APD}_{ref}|}.
\end{displaymath}

\section{Numerical results}
\label{sc:num_res}
\begin{table}[htbp!]%
  \centering
  \caption{Accuracy for the BR model for various classical and stabilized methods.}
  \label{rel_err}
  \subfigure[AB$_2$, RL$_2$, EAB$_2$ and CN ]{%
    \label{rel_err2}%
    \begin{tabular}{lllll}
      \toprule 
      $h$           &AB$_2$& RL$_2$      & EAB$_2$ & CN \\
      \midrule
      \num{0.2}     & -- & \num{0.251}   & \num{0.284}   & \num{4.11e-2} \\
      \num{0.1}     & -- & \num{0.107}   & \num{9.26e-2} & \num{1.13e-2} \\
      \num{0.05}    & -- & \num{3.35e-2} & \num{2.31e-2} & \num{2.65e-3} \\
      \num{0.025}   & -- & \num{8.88e-3} & \num{5.39e-3} & \num{6.66e-3} \\
      \num{0.0125}  & -- & \num{2.23e-3} & \num{1.29e-3} & \num{1.68e-4} \\
      \num{6.25e-3} & \num{2.07e-4} & \num{5.6e-4} & \num{3.17e-4} & \num{4.25e-5} \\
      \bottomrule
    \end{tabular}%
  }      
  \subfigure[AB$_3$, RL$_3$, EAB$_3$ and BDF$_3$ ]{%
    \label{rel_err3}%
    \begin{tabular}{lllll}
      \toprule
      $h$& AB$_3$ & RL$_3$ & EAB$_3$ & BDF$_3$ \\
      \midrule
      \num{0.2} & --  & \num{0.148} & 0.516        & \num{4.09e-2} \\
      \num{0.1} & --    & \num{4.07e-2} & \num{9.17e-2} & \num{1.04e-2} \\
      \num{0.05} & --     & \num{6.34e-3} & \num{1.09e-2} & \num{2.29e-3} \\
      \num{0.025} & --        & \num{7.57e-4} & \num{1.17e-3} & \num{3.84e-4} \\
      \num{0.0125} & --           & \num{9.07e-5} & \num{1.4e-4} & \num{5.25e-5} \\
      \num{6.25e-3} & \num{1.13e-5} & \num{8.23e-6} & \num{1.72e-5} & \num{2.01e-5} \\
      \bottomrule
    \end{tabular}%
  }    
  \subfigure[RK$_4$, RL$_4$, EAB$_4$ and BDF$_4$]{%
    \label{rel_err4} %
    \begin{tabular}{lllll}
      \toprule
      $h$       & RK$_4$   & RL$_4$       & EAB$_4$ & BDF$_4$ \\
      \midrule
      \num{0.2}    & --      & --            & --         & \num{4.98e-2} \\
      \num{0.1}    & --       & \num{5.86e-2}  & \num{0.119} & \num{1.27e-2} \\
      \num{0.05}    & --       & \num{4.58e-3}  & \num{8.96e-3} & \num{2.02e-3} \\
      \num{0.025} & \num{4.65e-5} & \num{2.61e-4} & \num{4.33e-4} & \num{1.93e-4} \\
      \num{0.0125} & \num{2.67e-6} & \num{1.62e-5} & \num{2.67e-5} & \num{3.52e-5} \\
      \num{6.25e-3} & \num{1.65e-7} & \num{9.94e-7} & \num{1.73e-6} & \num{2.01e-5} \\
      \bottomrule
    \end{tabular}%
  }
\end{table}
\subsection{Accuracy}
The relative error $e(h)$ is computed for various time-steps $h$ and depicted in Table \ref{rel_err} where it can be observed that all the methods exhibit the expected order of convergence. A general view of Table \ref{rel_err} shows that the RL$_k$ is always more accurate than EAB$_k$ and unlike the classical explicit schemes, the stabilized schemes allows the use of large time-steps as the implicit except at the order four where it is not possible for $h=0.2$.

Table \ref{rel_err}(a) shows that the CN is the most accurate among the methods of order 2 with a factor in the range of 10.
Table \ref{rel_err}(b) shows that the BDF$_3$ method is better than the stabilized schemes for $h\geq 0.0125$ with a coefficient 10 for $h=0.2$ while for $h<\num{6.25e-3}$ the RL$_3$ is more accurate. 
Table \ref{rel_err}(c) shows that the RK$_4$ method is the most accurate among the methods of order 4 for $h\leq 0.025$ while for $h>0.025$ the BDF$_4$ is more accurate than the stabilized schemes.  

\subsection{Cost}
\begin{figure}[htbp!]
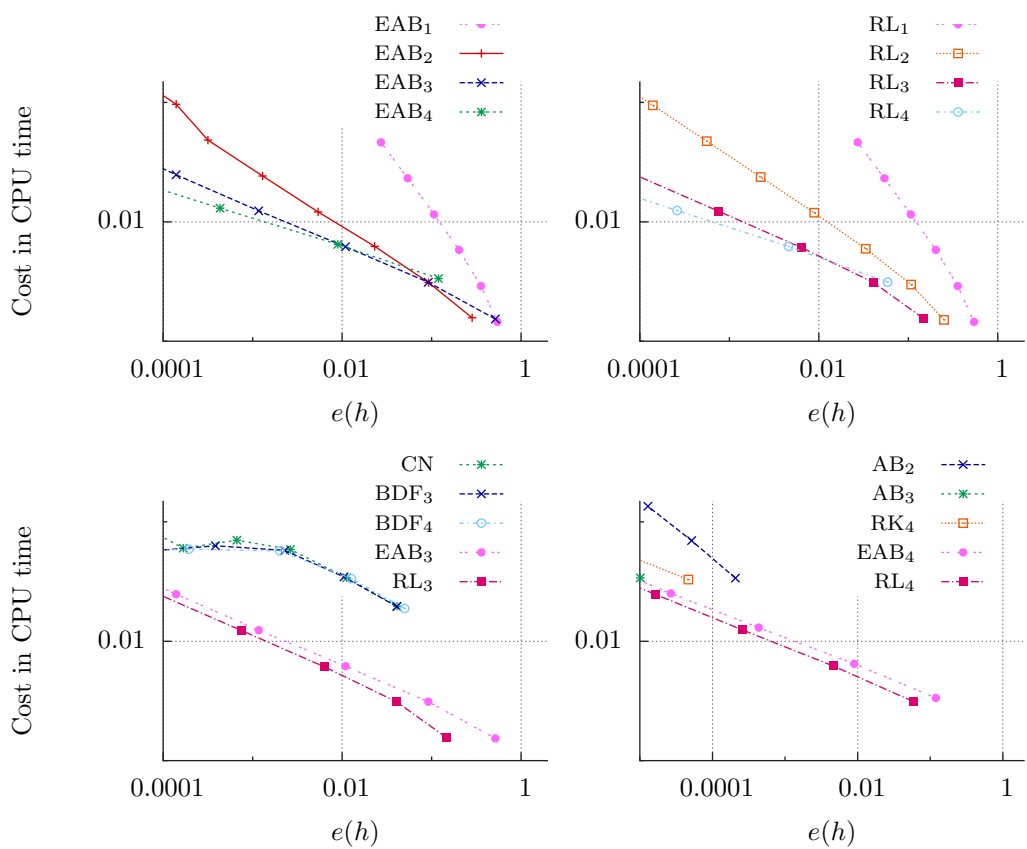

  \centering
  \input{EABk_cost.txt} \input{RLk_cost.txt} 
  \\[15pt]
  \input{RL3_EAB3_sch_imp_cost.txt} \input{RL4_EAB4_RKAB3_cost.txt}\\[5pt]
  \caption{CPU time as a function of the error $e(h)$ in Log/Log scale.}
  \label{fig:CPU-time}
\end{figure}

A general observation of Figure \ref{fig:CPU-time} on the top shows that for the error between $1\%$ and $10\%$ the gain in terms of CPU time is high (with a factor in the range of 10) when moving from the order 1 to order 2 schemes. This gain remains important (with a factor in the range of 5) when moving from the order 2 to order 3 schemes while for the errors between $1\%$ and $10\%$ there is no gain when moving from the order 3 to the order 4 schemes. However the order 4 becomes advantageous for the errors less than $0.1 \%$.

Figure \ref{fig:CPU-time} on the bottom shows that the RL$_3$ and the RL$_4$ are less costly than the EAB$_3$ and EAB$_4$ respectively. The factor is not so high but in terms of implementation, the RL is easier than the EAB schemes.

Figure \ref{fig:CPU-time} on the bottom left shows that when using high order stabilized schemes instead of implicit schemes, the gain in time CPU is very high with a coefficient greater than 10. This is due to the fact that the nonlinear solver is very expensive and its cost become very high for large time-steps.

Figure \ref{fig:CPU-time} on the bottom right shows that the order 4 stabilized schemes are less costly than the classical explicit schemes but it is better to use the RL$_4$ scheme instead of the EAB$_4$ scheme. Because of their stability properties the classical explicit schemes require the use of small time-steps that make them sometimes useless. For instance the RK$_4$ is very accurate but its use require to take a small time step. This small time steps produces a very small error that might be not needed and then its use will induce an additional cost.

\begin{figure}[htbp!]
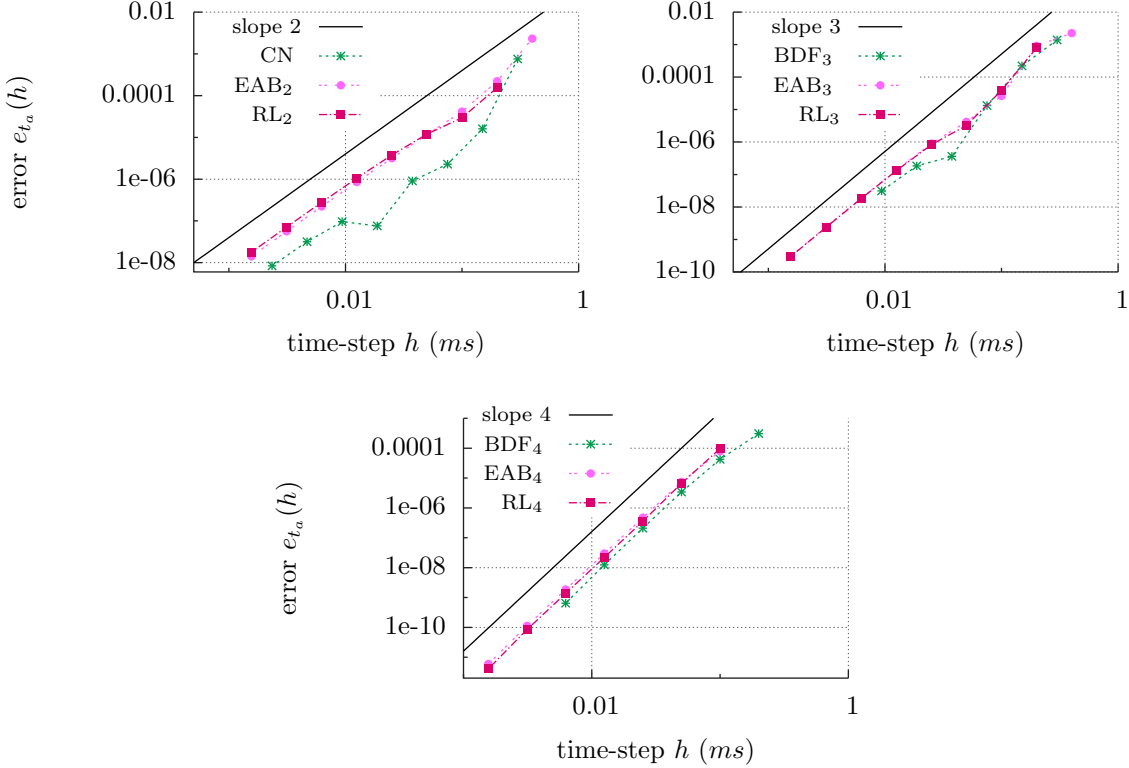

  \centering 
  \input{order2_ta_APD_conv.txt} ~~~~~~\input{order3_ta_APD_conv.txt} 
  \\[10pt]
  \input{order4_ta_APD_conv.txt} \\[10pt]
  \caption{Relative error $e_{t_a}(h)$ for the CN, BDF$_3$, BDF$_4$, EAB$_k$  and RL$_k$ schemes $k=2,3,4$  for the BR model.}
  \label{fig:ta_APD_conv}
\end{figure}
\subsection{Accuracy on $t_a$, $t_r$ and $APD$}
\begin{figure}[htbp!]
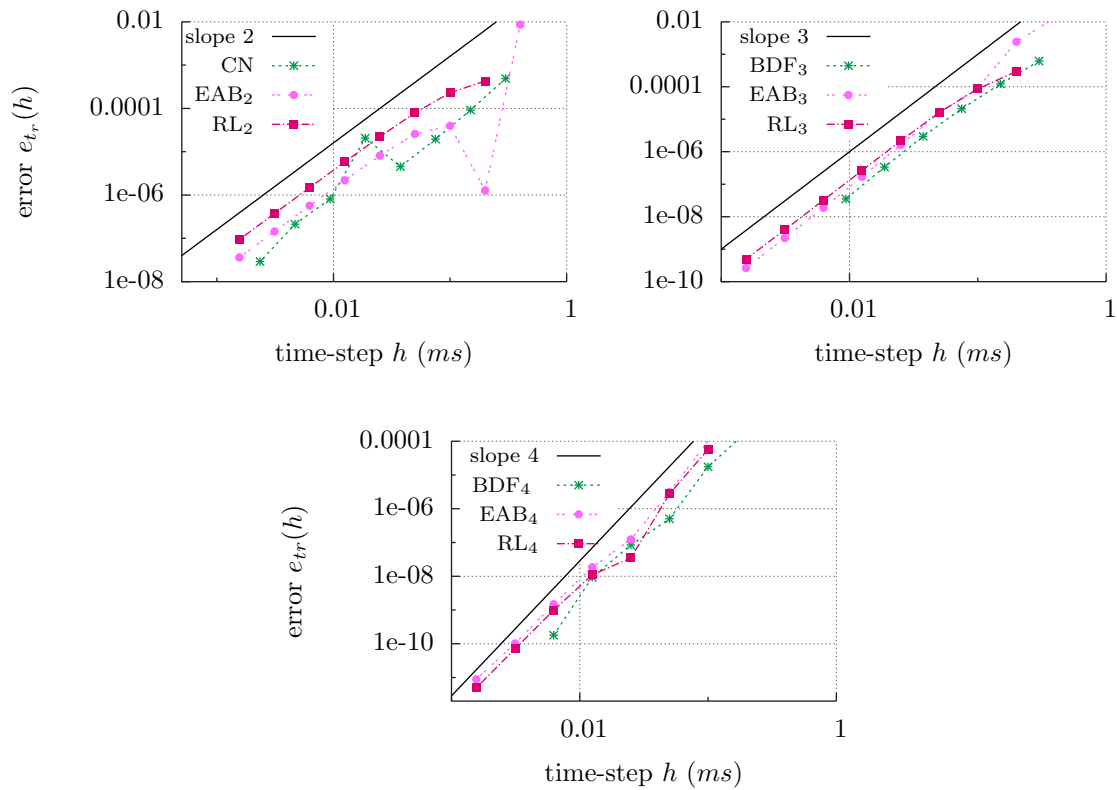

  \centering
  \input{order2_tr_APD_conv.txt} ~~~~~~\input{order3_tr_APD_conv.txt} 
  \\[15pt]
  \input{order4_tr_APD_conv.txt} \\[10pt]
  \caption{Relative error $e_{t_r}(h)$ for the CN, BDF$_3$, BDF$_4$, EAB$_k$  and RL$_k$ schemes $k=2,3,4$ for the BR model.}
  \label{fig:tr_APD_conv}
\end{figure}

We investigate in this section the accuracy on $t_a$, $t_r$ and $APD$. The previous results Section \ref{sc:num_res} showed that classical explicit methods induce very high computational costs because of their lack of stability. Only the implicit methods will then be considered to benchmark with the EAB$_k$ and RL$_k$ methods.
For a given numerical solution, the method described in Section \ref{sc:phys_quantities} is used to perform the values of $t_a$, $t_r$ and $APD$ . These values are then compared with the ones predicted by the reference solution. The errors for $t_a$ and $t_r$ are depicted on the figures \ref{fig:ta_APD_conv} and \ref{fig:tr_APD_conv} for various time-steps. These figures show that for a numerical solution computed with an order $k$ numerical scheme, the values of $t_a$ and $t_r$ predicted converge to the ones predicted by the reference solution with the same convergence order.
In the same figures, we can see that with equal time-step $h \leq 0.01$, the errors on the values predicted by the numerical solution decrease with a factor of $10$ at least, when moving from the order $k$ to the order $k+1$ schemes. 
We didn't show the pictures on $APD$ but since $APD = t_r - t_a$, the results on $APD$ will be the same as for $t_a$ and $t_r$.

\section{Conclusion}
\label{sc:conclusion}
Two families of explicit high order stabilized methods (EAB$_k$, RL$_k$) have been studied in this work. Excepted the order four, both have been shown to be as stable as the classical implicit methods for the test case we have chosen. Meanwhile the two families of schemes have been compared with some classical solvers (CN, BDF$2$, BDF$3$, AB$_2$, AB$_3$, RK$_4$). This comparison has shown (for the test case we chose)that EAB and RL are competitive when both the cost and the accuracy are taken in account.  Otherwise, it has also been demonstrated that the use of high order (3 or 4) of the stabilized methods instead of the classical high order implicit methods allows to decrease the cost almost 50 times. A method permitting to compute accurately (without degrading the convergence order of the numerical scheme) the values of $t_a$, $t_r$ and $APD$ predicted by a numerical solution has been also described and numerically investigated in this work.

\bibliographystyle{abbrv} \bibliography{biblio}

\end{document}